\documentclass[reqno]{amsart}
\usepackage{latexsym}
\usepackage{amssymb, color}
\theoremstyle{plain}
\usepackage{txfonts}
\usepackage[sc]{mathpazo}

\def\XXint#1#2#3{{\setbox0=\hbox{$#1{#2#3}{\int}$}
     \vcenter{\hbox{$#2#3$}}\kern-.5\wd0}}

\newcommand{\chara}{1\!\!1}

\newcommand{\na}{\nabla}

\newcommand{\lt}{\left}
\newcommand{\rt}{\right}
\newcommand{\nl}{\newline}

\newcommand{\nn}{\nonumber}

\newcommand{\lm}{\lambda}

\newcommand{\qd}{\quad}

\newcommand{\ep}{\epsilon}

\newcommand{\II}{\mathcal{I}}

\newcommand{\PPI}{\mathcal{P}}

\newcommand{\GI}{\mathcal{I}}

\newcommand{\UI}{\mathcal{u}}

\newcommand{\OI}{\mathcal{O}}

\newcommand{\ti}{\tilde}

\newcommand{\R}{\mathrm {I\!R}}

\newcommand{\dia}{\diamondsuit}

\newcommand{\red}[1]{{\textcolor{black}{#1}}}

\newtheorem{a1}{Lemma}
\newtheorem{a2}{Theorem}

\newtheorem{a6}{Corollary}

\newtheorem*{thm}{Theorem}

\theoremstyle{remark}

\begin{document}
\title[A generalized \red{Stoilow} decomposition for pairs of mappings of integrable dilatation]
{A generalized \red{Stoilow} decomposition for pairs of mappings of integrable dilatation}
\author{Andrew Lorent}
\address{Mathematics Department\\University of Cincinnati\\2600 Clifton Ave.\\ Cincinnati OH 45221 }
\email{lorentaw@uc.edu}
\maketitle

\begin{abstract}
We prove a rigidity result for pairs of mappings of integrable dilatation whose gradients pointwise deform the unit ball to similar ellipses. Our result implies as corollaries a version of the generalized 
\red{Stoilow} decomposition provided by Theorem 5.5.1 of \cite{astala} and the two dimensional rigidity result of \cite{lor20} 
for mappings whose symmetric part of gradient agrees. 

Specifically let $u,v\in W^{1,2}(\Omega,\R^2)$ where $\det(Du)>0$, $\det(Dv)>0$ a.e.\ and $u$ is a mapping of integrable dilatation. Suppose for 
a.e.\ $z\in \Omega$ we have \red{$Du(z)^T Du(z)=\lm Dv(z)^T Dv(z)$} for some $\lm>0$. Then there exists a meromorphic function $\psi$ and a homeomorphism $w\in W^{1,1}(\Omega:\R^2)$ such that 
$Du(z)=\PPI\lt(\psi(w(z))\rt)Dv(z)$ where $\PPI(a+ib)=\lt(\begin{smallmatrix} a & -b \\ b & a\end{smallmatrix}\rt)$.

We show by example that this result is sharp in the sense that there can be no continuous relation between the gradients of $Du$ and $Dv$ on a dense open \red{connected} subset of $\Omega$ unless one of the mappings is of integrable dilatation.   
\end{abstract}

In Theorem \red{5.5.1} of their recent monograph \cite{astala} Astala, Iwaniec, Martin provide a 
generalized \red{Stoilow} decomposition. Specifically they show that if $u\in W^{1,1}$, $v\in W^{1,2}$ are Quasiregular mappings, $u$ is a homeomorphism and the Beltrami coefficients of $u$, $v$ are the same, then $u=\phi\circ v$ 
for some holomorphic mapping $\phi$. This is 
a generalization of the classical \red{Stoilow} decomposition for quasiregular mappings\footnote{To see take $u$ to be a quasiregular mappings and $v$ to be the solution of the Beltrami equation with Beltrami coefficient of $u$ apply Theorem \red{5.5.1}}. In Theorem 20.2.1 \cite{astala} a still more general \red{Stoilow} theorem is proved that as a corollary has the generalization of Theorem \red{5.5.1} for mappings of integrable dilatation, where one of these mappings is a 
homeomorphism. A different generalization is provided by Theorem 20.4.19 where no assumption of integrability of 
the dilatation is needed, however slightly stronger integrability assumptions on $u$, $v$ are required.   

In \cite{cia}, \cite{cia1} a closely related problem has been studied from the context of elasticity. Cialet and Mardare considered deformations $u:\Omega\rightarrow \R^n$ whose determinant is positive everywhere and studied the relation between the `Cauchy Green' tensor $Du^T Du$ and the deformation. They proved a kind of continuity property of this relation that implies the well known result that if a pair of $C^1$ deformations have the same Cauchy Green tensor then they are related by an isometry, for a proof 
see \cite{cia2}, Theorem 3.

In a recent paper motivated by powerful results on rigidity of differential inclusions \cite{fmul} we proved a sharp generalization 
of the $C^1$ result characterizing $C^1$ functions whose Cauchy Green tensor agree. We showed that if $u,v\in W^{1,2} $ are mappings 
of integrable dilatation and the symmetric part of the gradient of $u$ and $v$ are the same, then $u=l_{R}\circ v$ with $Dl_{R}=R\in SO(2)$, see Theorem 1 \cite{lor20}. In this paper our main result will be a generalized \red{Stoilow} decomposition for pairs of $W^{1,1}$ mappings of positive 
determinant whose conformal part of gradient is square integrable where one of these mappings is of integrable dilatation. This result will imply Theorem 1 \cite{lor20} and the 
generalization of Theorem 5.5.1 \cite{astala} for mappings of integrable dilatation (defined on simply connected domains), as simple corollaries. As far as we are aware this is the only \red{Stoilow} decomposition result for pairs of mappings that does not assume one of 
them is injective. \red{Given matrix $A$, by polar decomposition we can write $A=R(A)S(A)$ where $R(A)$ is a rotation and 
$S(A)$ is a symmetric matrix given by $\sqrt{A^T A}$}. Our main theorem is:  
\begin{a2}
\label{T2}
Let $\Omega\subset \R^2$ be a connected open domain, let 
$v\in W^{1,1}(\Omega:\R^2)$ and \red{$u\in W^{1,1}(\Omega:\R^2)$}, $\det(D v)>0$, $\det(D u)>0$ a.e.\, $\int_{\Omega} \lt|\lt[Du\rt]_c\rt|^2+\lt|\lt[Dv\rt]_c\rt|^2 dx<\infty$ (where $\lt[\cdot \rt]_c$ denotes the conformal part of a matrix) and $\|D u(z)\|^2\leq K(z)\det(D u(z))$ for $K\in L^1$. Suppose: 
\begin{equation}
\label{bvbb51}
\text{ For a.e. }z\in \Omega,\;\;
S(Du(z))=\lm S(Dv(z)) \text{ for some }\lm>0.
\end{equation}
Then there exists a homeomorphism $w\in W^{1,1}(\Omega:\R^2)$ and meromorphic function $\psi:w(\Omega)\rightarrow \mathbb{C}$ such that 
\begin{equation}
\label{gheq99}
Dv(z)=\PPI\lt(\psi(w(z))\rt)Du(z)\text{ for }a.e.\ z\in \Omega
\end{equation}
where $\PPI$ is defined by $\PPI(a+ib)=\lt(\begin{smallmatrix} a & -b \\ b & a\end{smallmatrix}\rt)$. 
In addition the poles of $\psi$ are contained in $w(B_u)$ where $B_u$ denotes the branch set of $u$, 
(the set of points where 
$u$ fails to be a local homeomorphism).
\end{a2}

One of the principle corollaries is the two dimensional version of Theorem 1 \cite{lor20}.
\begin{a6}
\label{CC2}
Let $\Omega\subset \R^2$ be a connected open domain, let 
$v\in W^{1,1}(\Omega:\R^2)$ and $u\in W^{1,2}(\Omega:\R^2)$, $\det(Dv)>0$, $\det(Du)>0$ a.e.\ and \red{$\|Du\|^2\leq K\det(Du)$} for $K\in L^1$. Suppose 
\begin{equation}
\label{bvbb51.5}
Du^T Du=Dv^T Dv \text{ for }a.e.\ x\in \Omega
\end{equation}
then there exists \red{$R\in SO(2)$} 
\begin{equation}
\label{bvbb50}
Dv=R Du\text{ for }a.e.\ x\in \Omega.
\end{equation}
\end{a6}

We will show by Example 1 of Section \ref{example} that there can be no continuous relation between 
$Dv$ and $Du$ on a dense \red{connected} open subset of $\Omega$ unless one of $u$ or $v$ is of integrable dilatation. A slightly more general version of Example 1 was given in \cite{lor20} which showed the sharpness of the two dimensional version of Theorem 1 of \cite{lor20}. The example and Theorem 1 of \cite{lor20} 
also answered the question 
posed in \cite{cia} as to the the optimal function class for which 
(\ref{bvbb51.5}) implies (\ref{bvbb50}). We showed the optimal function class is the space of functions of integrable dilatation (note however an example was already presented in \cite{cia} which 
showed this implication could not hold for arbitrary functions in $W^{1,2}$). Theorem \ref{T2} of this paper and Example 1 of Section \ref{example} shows that the optimal class for the more general question of when (\ref{bvbb51}) implies a `global' continuous relation between $Du$ and $Dv$ is also the space of functions with integrable dilatation.

The following corollary is essentially folklaw, it follows from Theorem 20.2.1 \cite{astala}. We 
state it in matrix notation.
\begin{a6}
\label{CC4}
Let $\Omega\subset \R^2$ be a connected open domain, let 
$v\in W^{1,1}(\Omega:\R^2)$ and \red{$u\in W^{1,1}(\Omega:\R^2)$}, $\det(D v)>0$, $\det(D u)>0$ a.e.\, $\int_{\Omega} \lt|\lt[Du\rt]_c\rt|^2+\lt|\lt[Dv\rt]_c\rt|^2 dz<\infty$ (where $\lt[\cdot\rt]_c$ denotes the conformal part of the matrix) and $\|D u(z)\|^2\leq K(z)\det(D u(z))$ for $K\in L^1$. Suppose: 
\begin{equation}
\text{ For a.e. }z\in \Omega,\;\;
S(Du(z))=\lm S(Dv(z)) \text{ for some }\lm>0.\nn
\end{equation}
Then for any $z\in \Omega\backslash B_u$ (where $B_u$ denotes the branch set of $u$) there 
exists $r_z>0$ and conformal map $\phi:B_{r_z}(z)\rightarrow \R^2$ such that 
$$
v(x)=\phi(u(x))\text{ for all }x\in B_{r_z}(z). 
$$
\end{a6}

Another corollary is the following \red{Stoilow} decomposition for mappings of integrable dilatation which is a generalization of Theorem 5.5.1 \cite{astala} for simply connected domains but is a corollary to the more general result Theorem 20.2.1 \cite{astala}.

\begin{a6}
\label{CC3}
Let $\Omega\subset \mathbb{C}$ be a simply connected open domain, let 
$u\in W^{1,1}(\Omega:\mathbb{C})$ be a homeomorphism that satisfies 
$\frac{\partial u}{\partial z}\in L_{loc}^2(\Omega)$. Let $v\in W^{1,2}(\Omega:\mathbb{C})$. Suppose both $u$ and $v$ satisfy 
the Beltrami equation 
\begin{equation}
\frac{\partial u}{\partial \bar{z}}=\mu\frac{\partial u}{\partial z}\text{ and }\frac{\partial v}{\partial \bar{z}}=\mu\frac{\partial v}{\partial z}\text{ for }a.e.\ z\in \Omega
\end{equation}
where $\mu:\Omega\rightarrow \mathbb{C}$ is a measurable function such that $\lt|\mu(z)\rt|<1$ a.e.\ and 
$\int_{\Omega} \frac{1+\lt|\mu(z)\rt|}{1-\lt|\mu(z)\rt|} dz<\infty$ then there exists a holomorphic function $\phi:u(\Omega)\rightarrow \mathbb{C}$ such that  
\begin{equation}
v=\phi\circ u.\nn
\end{equation}
\end{a6}

\red{Note Corollary \ref{CC3} is a generalization of Theorem 5.5.1 of \cite{astala} (for simply connected domains) because it is a known result that 
a $W^{1,1}$ homeomorphism $u$ is such that $\det(Du)\in L_{loc}(\Omega)$ (for the convenience of the reader we 
include a proof of this result in the appendix). Thus as a quasiregular mapping will have the determinant of the 
conformal part of its gradient bounded by a constant multiple of its determinant so it follows that any quasiregular mapping 
$u\in W^{1,1}$ that is a homeomorphism will be such that $\lt[Du\rt]_c\in L^2_{loc}$ and hence satisfies the hypothesis of Corollary \ref{CC3}.}

Theorem \ref{T2} is best thought of a kind of rigidity result. Suppose we have two mappings \red{with pointwise} positive determinant  
$u$, $v$ and at a.e.\ point $z$ we have that the ellipses $E_{Du(z)}:=\lt\{Du(z)\phi:\phi\in B_1(0)\rt\}$ and 
$E_{Dv(z)}:=\lt\{Dv(z)\phi:\phi\in B_1(0)\rt\}$ are similar. Then for a.e.\ $z$ there exists a matrix $C(z)\in CO_{+}(2)$ such that 
$Du(z)=C(z)Dv(z)$. Theorem \ref{T2} says that if we know in addition that the 
conformal part of the gradient of the mappings are integrable and one of these mappings is of integrable dilatation 
then $C(z)$ is actually a continuous functions outside an isolated set of points. Indeed $C(z)=\PPI(\psi(w(z)))$ where 
$\PPI(a+ib)=\lt(\begin{smallmatrix} a & -b \\ b & a\end{smallmatrix}\rt)$, and $\psi$ is meromorphic and $w\in W^{1,1}$ is a 
homeomorphism. If neither of these mappings is of integrable dilatation then no such result is possible as we will show by a 
simple counterexample in Section \ref{example}. Even when both of these mappings is of integrable dilatation we can not expect $C(z)$ to be continuous everywhere. 
In complex notation, simply take $u(z)=z$ and $v(z)=(z-1)^2$ we find $Du(z)=C(z)Dv(z)$ where $\lim_{z\rightarrow 1} C=\infty$. 

 The most concise way in matrix notation to express the hypothesis that $E_{Du(z)}$ and $E_{Dv(z)}$ are similar ellipses is by equation (\ref{bvbb51}), i.e.\ 
to insist that the symmetric parts of $Du(z)$ and $Dv(z)$ are scalar multiples of each other. As mentioned in \cite{lor20} we considered 
the more specific question of the rigidity of pairs of mappings whose symmetric part of gradient agree. This question was partly  
motivated by the powerful rigidity and stability result for mappings whose gradient is close to $SO(n)$ proved in \cite{fmul}. It was 
shown in \cite{lor20} that if two mappings $u,v$ have pointwise positive determinant and one of them is of integrable dilatation 
and they satisfy $S(Du)=S(Dv)$ a.e.\ then $Du(z)=RDv(z)$ a.e. for some $R\in SO(n)$. This result holds in all dimension so 
long as the dilatation of one of these mappings is in $L^p$ for $p>n-1$ for $n\geq 3$ and the dilatation is in $L^1$ if $n=2$. The proof uses truncation theorems and the stability result 
of \cite{fmul} to bypass the lack of a chain rule. In this paper we achieve a much simpler proof of the two dimensional result by 
invoking the power of the `measurable Riemann mapping theorem' which allows us to solve the Beltrami equation $\frac{\partial u}{\partial \bar{z}}=\mu\frac{\partial u}{\partial z}$ for arbitrary measurable $\mu$ where $\lt|\mu\rt|\leq k<1$ a.e. 

The connection 
between Beltrami's equation and mappings that satisfy (\ref{bvbb51}) and (\ref{bvbb51.5}) is as follows. If we take a $\Omega\subset \mathbb{C}$ and a function $u:\Omega\rightarrow \mathbb{C}$ then define the $\R^2$ valued function $\ti{u}(x,y)=(\mathrm{Re}(u(x+iy)),\mathrm{Im}(u(x+iy)))$. Let $CO_{+}(2)$ denote the 
set of conformal $2\times 2$ matrices. And let $\lt[\cdot\rt]_M$ denote the homomorphism between $\mathbb{C}$ and $CO_{+}(2)$, so 
$\lt[a+ib\rt]_M=\lt(\begin{matrix} a & -b \\ b & a\end{matrix}\rt)$. Finally recall any $A\in M^{2\times 2}$ can be decomposed uniquely into a conformal and anticonformal matrix and denote these by $\lt[A\rt]_a$ and $\lt[A\rt]_c$ respectively. It is straight forward to see that $\lt[\frac{\partial u}{\partial z}\rt]_M=\lt[D\ti{u}\rt]_c$ and  
$\lt[\frac{\partial u}{\partial \bar{z}}\rt]_M \lt(\begin{matrix} 1 & 0 \\ 0 & -1\end{matrix}\rt)=\lt[D\ti{u}\rt]_a$. So the Beltrami coefficient $\mu(z)$ in 
the Beltrami equation $\frac{\partial u}{\partial \bar{z}}(z)=\mu(z)\frac{\partial u}{\partial z}(z)$ is 
the complex number (or in matrix notation a conformal matrix) that relates the anticonformal part of the gradient matrix to the conformal part of 
the gradient matrix. We can formulate this for any matrix $A\in M^{2\times 2}$ and so we define the Beltrami coefficient of the \em matrix \rm $A$
to be the $2\times 2$ conformal matrix $\mu_A$ that satisfies $\lt[A\rt]_a \GI=\mu_{A} \lt[A\rt]_c$ where 
$\GI:=\lt(\begin{matrix} 1 & 0 \\ 0 & -1\end{matrix}\rt)$. Note for any $\lambda>0$ it \red{is} immediate that $\mu_{\lambda A}=\mu_{A}$. It is a 
slightly longer calculation to see that for any $R\in SO(2)$ we have $\mu_{R A}=\mu_{A}$. So the Beltrami coefficient of the matrix $A$ 
does not detect dilatations or rotations. If we consider geometry of the ellipse $E_A:=\lt\{Av:v\in S^1\rt\}$ then the size of this ellipse 
is undetectable from the Beltrami coefficient $\mu_A$  and the direction in which the ellipse lies (the orientation of the ellipse) is 
undetectable from the Beltrami coefficient but the \em geometry \rm of the ellipse is  encoded in the relation between the 
anti-conformal and conformal parts of $A$ and so it entirely determined by $\mu_A$. Contrast this with the symmetric part of the matrix 
$S(A):=\sqrt{A^T A}$ which encodes both the geometry of $E_A$ and its size, but does not detect the orientation. So 
considering two mappings $u,v$ for which $S(Du)=S(Dv)$ is a much stronger hypothesis than considering two 
mappings for the Beltrami coefficients agree. Indeed given $A,B\in M^{2\times 2}$ if we only were interested in the \em geometry \rm of the ellipses $E_A$, $E_B$ we could postulate that $S(A)=\lambda S(B)$ for some $\lambda>0$. It turns out this 
is equivalent to $\mu_A=\mu_B$, (this is the content of Lemma \ref{LL1}).  So the hypothesis (\ref{bvbb51}) of our main theorem is equivalent to insisting that the Beltrami coefficient of $Du$ and $Dv$ agree which when written in complex notation is the 
statement that $u$, $v$ satisfy the same Beltrami equation. So our main theorem is also a \red{Stoilow} decomposition result for pairs of 
(non-invertible) mappings one of which is of integrable dilatation. As mentioned the counter example we exhibit in Section \ref{example} shows that no `global' \red{Stoilow} decomposition 
result is possible for mappings that do not have integrable dilatation. If one of the mappings is actually a homeomorphism 
stronger results are possible, see for example Theorem 20.4.19 \cite{astala}. 

We choose to express our statement and arguments in matrix language because it appears that the geometric content of the very efficient $\frac{\partial}{\partial z}$, $\frac{\partial}{\partial \bar{z}}$ notation is not widely known in the broader 
analysis and applied analysis community. Since the original results on stability and rigidity of differential inclusions \cite{fmul} 
were motived by elasticity and the initial investigations in \cite{lor20}, \cite{cia} were written in this context we 
prefer to write our arguments in a notation consistent with these works. However, as we will point out, many of our arguments are classical from the 
theory of mappings of integrable dilatation, specifically we use many ideas from the seminal paper \cite{iws1}. 

The value of this note consists of the following three things. Firstly we establish what is to our knowledge the first `global' \red{Stoilow} decomposition result 
for pairs of mappings where neither of them are assumed to be homeomorphisms. Secondly we show by example that mappings of integrable dilatation are the widest class of mappings for which such a result is possible\footnote{ Contrast this with Theorem 20.4.19 \cite{astala} where a \red{Stoilow} decomposition is proved for a pair of mappings where one of them is a homeomorphism without any assumptions of the integrability of the dilatation (however with slightly stronger integrability assumptions)}. Thirdly we demonstrate the very close connection between the line of generalization of 
\cite{fmul} started in \cite{lor20} and considered previously in the context of elasticity \cite{cia}, \cite{cia1} and the much studied topic of \red{Stoilow} decomposition.\nl

\bf Acknowledgements. \rm I would like to thank the referee for careful reading and very many excellent suggestions. The simplified Lemma \ref{LL10} and Lemma \ref{LL5} and the overall improvement of 
presentation of the paper are due to suggestions of the referee.\nl

\section{Conformal, Anti-conformal decomposition of $2\times 2$ matrices}

The following algebra identities are well known though typically stated in complex notation. Given $A\in M^{2\times 2}$ we can decompose $A$ uniquely into conformal and anti-conformal parts as 
follows 
\begin{equation}
\lt(\begin{matrix} a_{11} & a_{12} \\
a_{21} & a_{22} \end{matrix}\rt)=\frac{1}{2}\lt(\begin{matrix} a_{11}+a_{22} & -(a_{21}-a_{12}) \\
a_{21}-a_{12} & a_{11}+a_{22} \end{matrix}\rt)+\frac{1}{2}\lt(\begin{matrix} a_{11}-a_{22} & a_{21}+a_{12} \\
a_{21}+a_{12} & -(a_{11}-a_{22}) \end{matrix}\rt).\nn
\end{equation}
So for arbitrary matrix $A$ let 
$$
\lt[A\rt]_c:=\frac{1}{2}\lt(\begin{matrix} a_{11}+a_{22} & -(a_{21}-a_{12}) \\
a_{21}-a_{12} & a_{11}+a_{22} \end{matrix}\rt)\text{ and }
\lt[A\rt]_a:=\frac{1}{2}\lt(\begin{matrix} a_{11}-a_{22} & a_{21}+a_{12} \\
a_{21}+a_{12} & -(a_{11}-a_{22}) \end{matrix}\rt).
$$
It will often be convenient to write this decomposition as $A=\alpha R_{\theta}+\beta N_{\psi}$ 
where 
$$
R_{\theta}:= \lt(\begin{matrix} \cos{\theta} & -\sin{\theta} \\
\sin{\theta} & \cos{\theta} \end{matrix}\rt)\text{ and  }
N_{\psi}:= \lt(\begin{matrix} \cos{\psi} & \sin{\psi} \\
\sin{\psi} & -\cos{\psi} \end{matrix}\rt)
$$
 
Let $\GI:=\lt(\begin{smallmatrix} 1 & 0\\ 0 & -1\end{smallmatrix}\rt)$. 
The \bf Beltrami coefficient \rm of a matrix $A$ is a conformal matrix $\mu_A$ that relates the conformal and anti-conformal parts of $A$, it is defined by 
\begin{equation}
\label{eq2}
\lt[A\rt]_a \GI=\mu_{A} \lt[A\rt]_c.
\end{equation}

\subsection{Elementary algebraic properties}

\red{Given matrix $A\in M^{m\times n}$ let $\|A\|$ denote the operator norm of the matrix. Let $\lt|A\rt|$ denote the 
Hilbert Schmidt norm. Finally given complex number $z$ let $\lt|z\rt|$ denote the absolute value of $z$. } 
If $A=\alpha R_{\theta}+\beta N_{\psi}$ for $\alpha,\beta>0$ then 
\begin{equation}
\label{gheq50}
\det(A)=\det\lt(\lt[A\rt]_c\rt)+\det\lt(\lt[A\rt]_a\rt)=\alpha^2-\beta^2 
\end{equation}
and 
\begin{equation}
\label{bbreq942}
\det(\mu_A)=\frac{\beta^2}{\alpha^2},
\end{equation}
\begin{equation}
\label{brq1}
\|A\|\leq \alpha+\beta. 
\end{equation}
By (\ref{gheq50}), (\ref{bbreq942}) for any $A\in M^{2\times 2}$ with $\det(A)>0$ we have 
$0<\det(\mu_A)< 1$ and since $\mu_A$ is conformal 
\begin{equation}
\label{mmeq1}
0<\lt|\mu_A\rt|<\sqrt{2}\text{ for }A\in M^{2\times 2}, \det(A)>0.
\end{equation}
So 
\begin{equation}
\label{eq3}
\frac{\|A\|^2}{\det(A)}\overset{(\ref{brq1}),(\ref{gheq50})}{\leq} \frac{(\alpha+\beta)^2}{\alpha^2-\beta^2}=\frac{\alpha+\beta}{\alpha-\beta}. 
\end{equation}
Note that for any matrices $A,B\in M^{2\times 2}$ we have 
\begin{equation}
\label{breq400}
\lt[A+B\rt]_{a}=\lt[A\rt]_a+\lt[B\rt]_{a}\text{ and }\lt[A+B\rt]_{c}=\lt[A\rt]_c+\lt[B\rt]_{c}.
\end{equation}
It is immediate that if $C$ is a $2\times 2$ conformal matrix then $\lt[C\rt]_c=C$ and $\lt[C\rt]_a=0$ and 
if $B$ is a $2\times 2$ anticonformal matrix then $\lt[B\rt]_a=B$ and $\lt[B\rt]_c=0$. Thus for any matrix $A\in M^{2\times 2}$ and $C\in CO_{+}(2)$ we have 
\begin{equation}
\label{jkj1}
\lt[CA\rt]_c=\lt[C\lt[A\rt]_c+C\lt[A\rt]_a\rt]_c\overset{(\ref{breq400})}{=}\lt[C\lt[A\rt]_c\rt]_c=C \lt[A\rt]_c
\end{equation}
and 
\begin{equation}
\label{jkj2}
\lt[CA\rt]_a\overset{(\ref{breq400})}{=}\lt[C\lt[A\rt]_c+C\lt[A\rt]_a\rt]_a=C \lt[A\rt]_a.
\end{equation}
In the same way 
\begin{equation}
\lt[A C\rt]_c=\lt[A\rt]_c C\text{ and }\lt[A C\rt]_a=\lt[A\rt]_a C.\nn
\end{equation}
Now notice 
\begin{equation}
C\lt[A\rt]_a\GI\overset{(\ref{jkj2})}{=}\lt[C A\rt]_a\GI=\mu_{CA}\lt[CA\rt]_c\overset{(\ref{jkj1})}{=}C\mu_{CA}\lt[A\rt]_c\nn
\end{equation}
and so $\mu_{CA}=\mu_A$ and thus as previously mentioned, the Beltrami Coefficient does not detect dilatations or rotations.

Now letting $\alpha>0$, $\beta>0$, $\gamma,\lambda\in \lt[0,2\pi\rt)$ be such that $A=\alpha R_{\theta}+\beta N_{\psi}$ we have  
\begin{eqnarray}
\|A\|^2\leq Q\det{A}&\overset{(\ref{brq1}),(\ref{gheq50})}{\Rightarrow}&\frac{\lt(\alpha+\beta\rt)^2}{\alpha^2-\beta^2}
\leq Q\nn\\&\Rightarrow& (\alpha+\beta)\leq Q(\alpha-\beta)\nn\\&\Rightarrow& (1+\frac{\beta}{\alpha})\leq Q(1-\frac{\beta}{\alpha})\nn\\
&\Rightarrow&
\frac{\beta}{\alpha}
\leq\frac{Q-1}{Q+1}.\nn
\end{eqnarray}
As $\beta N_{\psi} \GI=\mu_{A} \alpha R_{\theta}$, so 
\begin{equation}
\lt|\mu_A\rt|=\frac{\beta}{\alpha}\nn
\end{equation}
hence  
\begin{equation}
\lt|\mu_A\rt|\rightarrow 1\text{ and 
}\frac{\|A\|^2}{\det(A)}\rightarrow \infty\text{ as }\beta\rightarrow \alpha.\nn
\end{equation}

\section{The relation between the Beltrami coefficient and symmetric part of a matrix}

\begin{a1} 
\label{L0.7}
Let $A\in M^{2\times 2}$, $\det(A)>0$. Let the Beltrami coefficient of $A$ be defined by (\ref{eq2}). The 
Beltrami coefficient of $A$ and $A^{-1}$ are related in the following way 
\begin{equation}
\label{eqw1.1}
\mu_{A} \lt[A\rt]_c \GI=-\mu_{A^{-1}} \GI \lt[A\rt]_c.
\end{equation}
\end{a1}

\em Proof of Lemma. \rm Let $\alpha>0$, $\beta>0$, $\theta,\psi\in \lt(0,2\pi\rt]$ be such that 
$A=\alpha R_{\theta}+\beta N_{\psi}$. Firstly note that 
\begin{equation}
\label{breq250}
R_{\theta} \GI= \lt(\begin{smallmatrix} \cos\theta & -\sin\theta \\ \sin\theta & \cos\theta \end{smallmatrix}\rt)
\lt(\begin{smallmatrix} 1 & 0 \\ 0 & -1 \end{smallmatrix}\rt)= \lt(\begin{smallmatrix} \cos\theta & \sin\theta \\ \sin\theta & -\cos\theta \end{smallmatrix}\rt)=\lt(\begin{smallmatrix} 1 & 0 \\ 0 & -1 \end{smallmatrix}\rt)\lt(\begin{smallmatrix} \cos\theta & \sin\theta \\ -\sin\theta & \cos\theta \end{smallmatrix}\rt)=\GI R_{-\theta}.
\end{equation}
We claim
\begin{equation}
\label{breq721}
A^{-1}=\frac{\alpha R_{-\theta}}{\alpha^2-\beta^2}-\frac{\beta N_{\psi}}{\alpha^2-\beta^2}.
\end{equation}
Note $N_{\psi}=R_{\psi}\GI$ so 
\begin{equation}
\label{kaq1}
N_{\psi} R_{-\theta}=R_{\psi} \GI R_{-\theta}\overset{(\ref{breq250})}{=}R_{\psi} R_{\theta} \GI=R_{\theta}N_{\psi}
\end{equation}
hence
\begin{eqnarray}
A A^{-1}&=& \lt(\alpha R_{\theta}+\beta N_{\psi}\rt)\frac{\lt(\alpha R_{-\theta}-\beta N_{\psi}\rt)}{\alpha^2-\beta^2}\nn\\
&=& \lt(\alpha^2-\beta^2\rt)^{-1}\lt( \alpha^2 Id+\alpha \beta N_{\psi}R_{-\theta}-\alpha \beta R_{\theta}N_{\psi}-\beta^2 N_{\psi}N_{\psi}\rt)\nn\\
&\overset{(\ref{kaq1})}{=}&\lt(\alpha^2-\beta^2\rt)^{-1}\lt( \alpha^2 Id-\beta^2 R_{\psi}\GI R_{\psi}\GI\rt)\nn\\
&\overset{(\ref{breq250})}{=}&\lt(\alpha^2-\beta^2\rt)^{-1}\lt( \alpha^2 Id-\beta^2 Id\rt) \nn\\
&=&Id. \nn 
\end{eqnarray}
In the same way we can see $A^{-1}A=Id$ and so (\ref{breq721}) is established. 

So as the decomposition into conformal and anticonformal parts in unique, by definition (\ref{eq2}) and from (\ref{breq721}) we have 
$\frac{-\beta N_{\psi}\GI}{\alpha^2-\beta^2}=\mu_{A^{-1}} \frac{\alpha R_{-\theta}}{\alpha^2-\beta^2}$ 
thus $N_{\psi}=-\mu_{A^{-1}}\frac{\alpha}{\beta}R_{-\theta}\GI$ and from definition (\ref{eq2}) 
$N_{\psi}=\mu_{A}\frac{\alpha}{\beta}R_{\theta}\GI$ so together these equations imply $\mu_{A} R_{\theta}=-\mu_{A^{-1}} R_{-\theta}$ and so 
\begin{equation}
\label{eq5.6}
\mu_{A}=-\mu_{A^{-1}}R_{-2\theta}.
\end{equation}
Hence $\|\mu_A\|=\|\mu_{A^{-1}}\|$ and so as $\mu_A$ and $\mu_{A^{-1}}$ are conformal 
\begin{equation}
\label{eq5.2}
\lt|\mu_{A^{-1}}\rt|=\lt|\mu_A\rt|.
\end{equation}
Now 
$$
\mu_{A}R_{\theta}\GI\overset{(\ref{eq5.6})}{=}-\mu_{A^{-1}}R_{-2\theta}R_{\theta}\GI=-\mu_{A^{-1}}R_{-\theta}\GI \overset{(\ref{breq250})}{=}   -\mu_{A^{-1}}\GI R_{\theta}
$$ 
which gives (\ref{eqw1.1}). $\Box$

%
%

\begin{a1}
\label{LL1}
Let $A,B\in M^{2\times 2}$ be matrices of positive determinant, then  
\begin{equation}
\label{fleq1}
\mu_A=\mu_B\Leftrightarrow S(A)=\lm S(B)\text{ for some }\lm>0.
\end{equation}
\end{a1}
\em Proof of Lemma \ref{LL1}. \rm Note the following equalities 
\begin{eqnarray}
\label{bbreq941}
\lt[A B^{-1}\rt]_a&:=&\lt[\lt(\lt[A\rt]_c+\lt[A\rt]_a\rt)\lt(\lt[B^{-1}\rt]_c+\lt[B^{-1}\rt]_a\rt)\rt]_a\nn\\
&=&\lt[A\rt]_a \lt[B^{-1}\rt]_c+\lt[A\rt]_c \lt[B^{-1}\rt]_a\nn\\
&\overset{(\ref{eq2})}{=}&\mu_A \lt[A\rt]_c \II \lt[B^{-1}\rt]_c+\lt[A\rt]_c \mu_{B^{-1}} \lt[B^{-1}\rt]_c\II\nn\\
&\overset{(\ref{eqw1.1})}{=}&\mu_A \lt[A\rt]_c \II \lt[B^{-1}\rt]_c-\lt[A\rt]_c \mu_{B} \II \lt[B^{-1}\rt]_c\nn\\
&=&(\mu_A-\mu_B)\lt[A\rt]_c \II \lt[B^{-1}\rt]_c.
\end{eqnarray}
So if $\mu_A=\mu_B$ then by (\ref{bbreq941}) $\lt[A B^{-1}\rt]_a=0$ so $A B^{-1}\in CO_{+}(2)$, i.e.\ 
there exists $k>0$ and $R\in SO(2)$ such that $A B^{-1}=k R$. So 
$A^T A=k^2 B^T B$ which implies $S(A)=k S(B)$. 

\red{Now on} the other hand if $S(A)=\lm S(B)$ for some $\lm>0$ then \red{recall by polar decomposition 
we have $A=R(A)S(A)$, $B=R(B)S(B)$ where $R(A),R(B)\in SO(2)$ and so} 
$$
A B^{-1}=R(A)S(A)\lt(R(B)S(B)\rt)^{-1}=\lm^{-1} R(A)(R(B))^{-1}\in CO_{+}(2). 
$$
So $\lt[A B^{-1}\rt]_a=0$ and so by (\ref{bbreq941}) we have that $\mu_A=\mu_B$ which concludes 
the proof of (\ref{fleq1}). $\Box$

\section{The Beltrami equation}
\label{beltrami}

As described in the introduction the Beltrami equation is a linear complex PDE \red{that} relates the conformal part of the gradient to the 
anti-conformal. So given a function $f$ from the complex plane to itself, $f(x+iy)=u(x,y)+iv(x,y)$. 
As is standard, $\frac{\partial f}{\partial \bar{z}}=\frac{1}{2}(\partial_x+i\partial_y) f$ and 
$\frac{\partial f}{\partial z}=\frac{1}{2}(\partial_x-i\partial_y) f$.

The basic theorem about the solvability of the Beltrami equation (sometimes known as the measurable Riemann mapping theorem) \red{is the following theorem, see Theorem 5.1.1, Theorem 5.3.2 \cite{astala} and 
for the original papers \cite{boj1}, \cite{morrey2}} 
\begin{thm}[Bojarsky, Morrey]
\label{Q1}
Suppose that $0\leq \kappa<1$ and that $\lt|\mu(z)\rt|\leq \kappa\chara_{B_r}(z)$, $z\in\mathbb{C}$. Then 
there is a unique $f\in W^{1,p}_{loc}(\mathbb{C})$  for every $p\in\lt[2,1+\frac{1}{\kappa}\rt)$ such that 
\begin{eqnarray}
\label{sstt1}
&~&\frac{\partial f}{\partial \overline{z}}=\mu(z)\frac{\partial f}{\partial z}\text{ for almost every }z\in\mathbb{C}\nn\\
&~&f(z)=z+O(\frac{1}{z})\text{ as }z\rightarrow \infty.
\end{eqnarray}
\red{And $f$ is a $\frac{1+k}{1-k}$-quasiconformal homeomorphism of $\mathbb{C}$}. 
\end{thm}

This is known as the \em principle solution \rm of the Beltrami equation. As described in the introduction, we will 
use the following notation. Given complex valued function $f$ of a complex variable 
let 
\begin{equation}
\ti{f}(x,y):=\lt(\mathrm{Re}(f(x+iy)),\mathrm{Im}(f(x+iy))\rt).\nn
\end{equation}
Rewriting the Beltrami equation 
in matrix notation we have that $\ti{f}$ satisfies 
\begin{equation}
\lt[D\ti{f}\rt]_a\GI=\mu_{D\ti{f}}\lt[D\ti{f}\rt]_{c}.\nn
\end{equation}

\section{Proof of Theorem \ref{T2}}

\begin{a1}
\label{LL4}
Let $\Omega\subset \R^2$ be a connected open set and let 
$v\in W^{1,1}(\Omega:\R^2)$ and \red{$u\in W^{1,1}(\Omega:\R^2)$}, $\det(D v)>0$, $\det(D u)>0$ a.e.\, $\int_{\Omega} \lt|\lt[Du\rt]_c\rt|^2+\lt|\lt[Dv\rt]_c\rt|^2 dx<\infty$ and $\|Du\|^2\leq K\det(Du)$ for $K\in L^1$. Suppose: 
\begin{equation}
\label{bvbb51.4}
\text{ For a.e. }z\in \Omega,\;\;S(Du(z))=\lm S(Dv(z)) \text{ for some }\lm>0
\end{equation}
then there exists homeomorphism $w:\Omega\rightarrow \R^2$ where $w^{-1}\in W^{1,2}$ and holomorphic functions 
$\phi_u$, $\phi_v$ such that 
\begin{equation}
\label{obbeq30}
u=\phi_u\circ w\text{ and }v=\phi_v\circ w.
\end{equation}
\end{a1}

\em Proof of Lemma \ref{LL4}. \rm Though we have slightly weaker assumptions, most of the proof of this lemma comes from following very closely the proof of Theorem 1 \cite{iws1}. 
The essential point is that the \red{Stoilow} decompositions of $u,v$ provided by Theorem 1 of \cite{iws1} have the same homeomorphism. This fact would 
essentially be immediate for those familiar with the methods of \cite{iws1} (the homeomorphism is the inverse of the limit of solutions to the 
Beltrami equation whose Beltrami coefficient is a truncation of the Beltrami coefficient of $\mu_{Du}$, and of course $\mu_{Du}=\mu_{Dv}$). However for completeness we provide the details. Now 
\begin{equation}
\label{jeq2}
\frac{\|A\|^2}{\det(A)}\overset{(\ref{eq3})}{=}\frac{1+\frac{\beta}{\alpha}}{1-\frac{\beta}{\alpha}}
\overset{(\ref{bbreq942})}{=}\frac{1+\sqrt{\det(\mu_A)}}{1-\sqrt{\det(\mu_A)}}.\nn
\end{equation}
Now by (\ref{bvbb51.4}) and (\ref{fleq1}) of Lemma \ref{LL1} we know 
\begin{equation}
\label{jeq5}
\mu_{Du(z)}=\mu_{Dv(z)}.
\end{equation}
So 
\begin{eqnarray}
\infty&>&\int_{\Omega} \frac{\|Du\|^2}{\det(Du)} dz=\int_{\Omega} \lt(\frac{1+\sqrt{\det(\mu_{Du(z)})}}{1-\sqrt{\det(\mu_{Du(z)})}}\rt) dz\nn\\
&\overset{(\ref{jeq5})}{=}&\int_{\Omega} \lt(\frac{1+\sqrt{\det(\mu_{Dv(z)})}}{1-\sqrt{\det(\mu_{Dv(z)})}}\rt) dz=\int_{\Omega}  \frac{\|Dv\|^2}{\det(Dv)} dz.\nn
\end{eqnarray}
Thus $v$ is a mapping of integrable dilatation. 

Now \red{ note we are carrying out these arguments in matrix notation and the `matrix' Beltrami coefficient 
is a $2\times 2$ conformal matrix of (Hilbert Schmidt) norm less than $\sqrt{2}$ (recall (\ref{mmeq1})). So  
we will define new Beltrami coefficients that are truncated on the points $z$ where $\lt|\mu_{Du(z)}\rt|$ and (respectively) 
$\lt|\mu_{Dv(z)}\rt|$ are close to $\sqrt{2}$}
\begin{equation}
\label{obbeq1}
\mu_u^{\ep}(z):=\lt\{\begin{array}{lcl}  \mu_{Du(z)} & \text{ if }&\lt|\mu_{Du(z)}\rt|\leq \sqrt{2}-\ep, \\ 
   (\sqrt{2}-\ep) \frac{\mu_{Du(z)}}{\lt|\mu_{Du(z)}\rt|} & \text{ if }&\lt|\mu_{Du(z)}\rt|> \sqrt{2}-\ep,\\
0 & \text{ if }&z\not \in \Omega \end{array}\rt.\ 
\end{equation}
and 
$$
\mu_v^{\ep}(z):=\lt\{\begin{array}{lcl}  \mu_{Dv(z)} & \text{ if }&\lt|\mu_{Dv(z)}\rt|\leq \sqrt{2}-\ep, \\ 
   (\sqrt{2}-\ep) \frac{\mu_{Dv(z)}}{\lt|\mu_{Dv(z)}\rt|} & \text{ if }&\lt|\mu_{Dv(z)}\rt|> \sqrt{2}-\ep,\\
0 & \text{ if }&z\not \in \Omega.
 \end{array}\rt.\
$$
Now let $w_u^{\ep}$ be the principle solution of the Beltrami equation with Beltrami coefficient $\mu_u^{\ep}$ and let $w_v^{\ep}$ be the principle solution of the Beltrami equation with Beltrami coefficient $\mu_v^{\ep}$. \nl

\em Step 1. \rm We will show $w_u^{\ep}=w_v^{\ep}$ a.e.\ and letting $h^{\ep}=(w_u^{\ep})^{-1}$ 
for some subsequence $\ep_k\rightarrow 0$ we have 
\begin{equation}
h^{\ep_k}\overset{W^{1,2}_{loc}}{\rightharpoonup} h\text{ as }k\rightarrow \infty\nn
\end{equation}
and $h$ is an monotone and $\Omega'=h^{-1}(\Omega)$ is an open connected set. 

\em Proof of Step 1. \rm Firstly by Lemma \ref{LL1} we have that $\mu_u^{\ep}=\mu_v^{\ep}$ a.e.\ so we have that 
the principle solutions defined by these Beltrami coefficients are the same, i.e. $w^{\ep}_u=w^{\ep}_v$ and so 
as we have defined $h^{\ep}$ to be the inverse of $w^{\ep}_u$ it is also the inverse of $w^{\ep}_v$. 

Now by Proposition 2.2 \cite{iws1} (or by a brief calculation) we know that 
\begin{equation}
\label{obbeq15}
\int_{E} \lt|Dh^{\ep}\rt|^2 dx\leq \int_{h^{\ep}(E)} \frac{\|Du\|^2}{\det(Du)} dx.
\end{equation} 
Recall $h^{\ep}$ is the principle solution to the Beltrami equation, \red{so note from the statement of Theorem \ref{Q1} we know that if 
$z\not\in \Omega$ then $\frac{\partial h^{\ep}}{\partial \bar{z}}(z)=0$. Since from (\ref{sstt1}) we 
know $h^{\ep}(z)=z+O(\frac{1}{z})$ so defining $\ti{h}^{\ep}(z)=h^{\ep}(z)-z$ we have that $\tilde{h}^{\ep}$ is 
analytic in $\mathbb{C}\backslash B_{\mathrm{diam}(\Omega)}(0)$. Finally define $g^{\ep}(z)=\tilde{h}^{\ep}(\frac{1}{z})$, so $g^{\ep}$ is analytic in $B_{1/\mathrm{diam}(\Omega)}(0)$ and 
$\lt|\lim_{z\rightarrow 0} g^{\ep}(z)z\rt|=\lt|\lim_{z\rightarrow 0} \tilde{h}^{\ep}(\frac{1}{z})z\rt|=0$ so 
$g^{\ep}$ has a removable singularity at $0$ and thus has a Talyor expansion around $0$. 
So $g^{\ep}(z)=\sum_{k=1}^{\infty} a_k z^k$ for $z\in B_{1/\mathrm{diam}(\Omega)}(0)$. Hence  
$h^{\ep}(\frac{1}{z})=\frac{1}{z}+\sum_{k=1}^{\infty} a_k z^k$ so changing variables, $h^{\ep}$ has a 
Talyor expansion at $\infty$} given by $h^{\ep}(z)=z+a_1 z^{-1}+a_2 z^{-2}+\dots$. Since $h^{\ep}$ is 
injective by Koebe distortion inequality applied to $\lt(h^{\ep}\lt(1/z\rt)\rt)^{-1}$ it can be shown that 
\begin{equation}
h^{\ep}(B_{r+1}(0))\subset B_{r+3}(0)\text{ for every }r>0\nn
\end{equation}
and so (\ref{obbeq15}) implies $Dh^{\ep}\in L^2_{loc}$ with uniform bounds independent of $\ep$. 

Now by Proposition 2.2 \cite{iws1} we have that $h^{\ep}$ is equicontinuous on every ball $B_r$ 
with modulus of continuity $\sqrt{\log\lt(\frac{4}{\lt|z_1-z_1\rt|}\rt)}$ (for $\lt|z_1-z_2\rt|<2$) with \red{bound} depending only 
on $\int_{B_{r+1}} \frac{\|Du\|^2}{\det(Du)} dz$. So specifically the bound in independent of $\ep$ and we 
have an equicontinuous sequence. Indeed Proposition 2.2.\ \cite{iws1} also establishes that $d(h_{\ep}(z),\infty)<\frac{10}{1+\lt|z\rt|}$ (where $d$ is 
the cord arc metric) and so the sequence $h_{\ep}$ is equicontinuous with respect to the chord arc metric on the extended complex plane 
$\hat{\mathbb{C}}$. So by Ascoli Arzela we must be able to extract a subsequence $h^{\ep_k}$ that converges uniformly to $h$ with 
respect to the chord arc metric on $\hat{\mathbb{C}}$.

Since each $h^{\ep}$ is a homeomorphism it has the property that $\lt(h^{\ep}\rt)^{-1}(x)$ is connected for every $x$ which is to say $h^{\ep}$ is monotone. As noted in \cite{iws1} by a result of 
Kuratowski Lacher and Whyburn \cite{mcau} the set of monotone maps is a closed subset of the set of continuous 
functions from $X$ to $Y$ under uniform convergence so long as $Y$ is locally connected. Hence this implies that the limiting map $h$ we obtain from subsequence $h^{\ep_k}$ will also be monotone. Since $h$ is continuous this implies $\Omega'=h^{-1}(\Omega)$ is an open connected set. \nl

\em Step 2. \rm We will show $\phi_u:=u\circ h$ and $\phi_v:=v\circ h$ are holomorphic. 

\em Proof of Step 2. \rm The key point is that on the set \red{$G_{\ep}:=\lt\{z:\lt|\mu_{Du(z)}\rt|\leq \sqrt{2}-\ep\rt\}$} the 
Beltrami coefficient of \red{$w^{\ep}_u$} and $u$ are the same so by Lemma \ref{LL1} we know for any $z\in G_{\ep}$ there exists $\lm_z>0$ such that 
$S(Du(z))=\red{\lm_z S(Dw^{\ep}_u(z))}$. Hence for any \red{$y\in w_u^{\ep}(G_{\ep})$} we have 
$$
D\phi^{\ep}_u(y)=Du(h^{\ep}(y))Dh^{\ep}(y)=Du(h^{\ep}(y))\red{Dw^{\ep}_u(h^{\ep}(y))^{-1}}\in CO_{+}(2).
$$
As $G_{\ep_1}\subset G_{\ep_2}$ for some $\ep_1<\ep_2$ and $\lt|\mathbb{C}\backslash \lt(\bigcup_{\ep>0} G_{\ep}\rt)\rt|=0$. So 
we know that on increasingly large sets $D\phi^{\ep}_u$ is conformal. However to actually prove the limiting map is conformal we 
need an estimate of the form
\begin{equation}
\label{obbeq51}
\int_{\Pi} \mathrm{dist}^2\lt(D\phi^{\ep}_u,CO_{+}(2)\rt) dz\rightarrow 0\text{ as }\ep\rightarrow 0 
\end{equation}
for any $\Pi\subset \subset \Omega'$. So note as $D\phi^{\ep}_u=\lt[D\phi^{\ep}_u\rt]_a+\lt[D\phi^{\ep}_u\rt]_c$ we know $\mathrm{dist}(D\phi^{\ep}_u,CO_{+}(2))\leq \lt|\lt[D\phi^{\ep}_u\rt]_a\rt|$. Thus \red{writing (\ref{obbeq51})} as a complex linear equation what we need is 
\begin{equation}
\label{obbeq52}
\frac{\partial \phi_u^{\ep}}{\partial \bar{z}}=h_{\ep}\text{ where }\|h_{\ep}\|_{L^2(\Pi)}\rightarrow 0\text{ as }\ep\rightarrow 0.
\end{equation}
Now if we could establish an upper bound of the form 
\begin{equation}
\label{obbeq53}
\sup_{\ep>0} \int_{\Pi} \lt|D\phi^{\ep}_u\rt|^2 dz<C
\end{equation}
we would know that for some $\ep_{k}\rightarrow 0$ the subsequence $\phi_u^{\ep_k}$ will converge weakly in $W^{1,2}_{loc}(\Omega)$ 
to $\phi_u$ and 
linearity of (\ref{obbeq52}) will ensure that $\frac{\partial \phi_u^{\ep}}{\partial \bar{z}}=0$ weakly and so by Weyl's lemma (see for example Lemma A.6.10 \cite{astala}) the limiting 
map is holomorphic. 

So what is required is to establish (\ref{obbeq51}), (\ref{obbeq53}) and this can be done with calculations involving 
(\ref{eqw1.1}) of Lemma \ref{L0.7} as is indicated in \cite{iws1}. For completeness we include these calculations. 

So note 
\begin{eqnarray}
\label{mmeq3}
&~&\red{\lt[D\phi_u^{\ep}(z)\rt]_a}\nn\\
&~&\qd=\lt[Du(h^{\ep}(z))\rt]_a \lt[D w^{\ep}_u(h^{\ep}(z))^{-1}\rt]_c+
\lt[Du(h^{\ep}(z))\rt]_c \lt[D w^{\ep}_u(h^{\ep}(z))^{-1}\rt]_a\nn\\
&~&\qd\overset{(\ref{eq2})}{=}\mu_{Du(h(z))} \lt[Du(h^{\ep}(z))\rt]_c \GI \lt[D w^{\ep}_u(h^{\ep}(z))^{-1}\rt]_c
+\lt[Du(h^{\ep}(z))\rt]_c \mu_{Dw^{\ep}_u(h^{\ep}(z))^{-1}}\lt[D w^{\ep}_u(h^{\ep}(z))^{-1}\rt]_c \GI\nn\\
&~&\qd\overset{(\ref{eqw1.1})}{=}\lt(\mu_{Du(h(z))}-\mu_{D w^{\ep}_u(h(z))}\rt) \lt[Du(h^{\ep}(z))\rt]_c \GI \lt[D h^{\ep}(z)\rt]_c.
\end{eqnarray}
Now suppose for $y\in \Omega$ is such that $\lt|\mu_{Du(y)}\rt|>\sqrt{2}-\ep$ we have 
\begin{eqnarray}
\label{mmeq12}
\det\lt(\mu_{Du(y)}-\mu_{Dw^{\ep}_u(y)}\rt)&=&\det\lt(\lt(1- \lt(\frac{\sqrt{2}-\ep}{\lt|\mu_{Du(y)}\rt|}\rt)\rt)\mu_{Du(y)}\rt)\nn\\
&\leq&\lt(1-\frac{\sqrt{2}-\ep}{\lt|\mu_{Du(y)}\rt|}\rt)^2\det(\mu_{Du(y)})\nn\\
&\overset{(\ref{mmeq1})}{\leq}& c\ep^2.
\end{eqnarray}
Let $\Pi\subset \subset \Omega'$. Hence 
\begin{eqnarray}
\label{gheq42.3}
\int_{\Pi} \lt|\det\lt(\lt[D\phi_u^{\ep}\rt]_a\rt)\rt| dz&\overset{(\ref{mmeq12}),(\ref{mmeq3})}{\leq}& 
c\ep^2 \int_{\Pi} \det\lt(\lt[Du(h^{\ep}(z))\rt]_c\rt)\det\lt(\lt[Dh^{\ep}(z)\rt]_c\rt) dz.
\end{eqnarray}
Now $Dh^{\ep}(z)=\lt[D w_u^{\ep}(h^{\ep}(z))\rt]^{-1}$ so 
\begin{equation}
\lt|\mu_{Dh^{\ep}(z)}\rt|=\lt|\mu_{\lt(D w_u^{\ep}(h^{\ep}(z))\rt)^{-1}}\rt|\overset{(\ref{eq5.2})}{=}
\lt|\mu_{\lt(D w_u^{\ep}(h^{\ep}(z))\rt)}\rt|\overset{(\ref{obbeq1})}{\leq}\sqrt{2}-\ep.\nn
\end{equation}
\red{Note that for a conformal matrix $C=\lt(\begin{smallmatrix} a & -b\\ b & a \end{smallmatrix}\rt)$ 
we have  $\mathrm{det}(C)=a^2+b^2=2\lt|C\rt|^2$ so using this for the first inequality}
\begin{eqnarray}
\label{gheq40}
\det(\mu_{Dh^{\ep}(z)})&=&\frac{1}{2}\lt|\mu_{Dh^{\ep}(z)}\rt|^2\leq \frac{1}{2}(2-2\sqrt{2}\ep+\ep^2)\nn\\
&\leq& 1-\sqrt{2}\ep+\frac{\ep^2}{2}\leq 1-\ep.
\end{eqnarray}
Thus
\begin{eqnarray}
\label{gheq41}
-\det\lt(\lt[Dh^{\ep}(z)\rt]_a\rt)&\overset{(\ref{eq2})}{\leq}&\det\lt(\mu_{D h^{\ep}(z)}\rt)\det\lt(\lt[Dh^{\ep}(z)\rt]_c\rt)\nn\\
&\overset{(\ref{gheq40})}{\leq}& (1-\ep)\det\lt(\lt[Dh^{\ep}(z)\rt]_c\rt).
\end{eqnarray}
Hence 
\begin{eqnarray}
\label{gheq42}
\det\lt(Dh^{\ep}(z)\rt)&\overset{(\ref{gheq50})}{=}&\det\lt(\lt[Dh^{\ep}(z)\rt]_c\rt)+ \det\lt(\lt[Dh^{\ep}(z)\rt]_a\rt)\nn\\
&\overset{(\ref{gheq41})}{\geq}& \ep \det\lt(\lt[Dh^{\ep}(z)\rt]_c\rt).
\end{eqnarray}
Now 
\begin{eqnarray}
\label{gheq52}
\int_{\Pi} \lt|\det\lt(\lt[D\phi_u^{\ep}(z)\rt]_a\rt)\rt| dz&\overset{(\ref{gheq42.3}),(\ref{gheq42}) }{\leq}&
\ep \int_{\Pi} \det\lt(\lt[Du(h^{\ep}(z))\rt]_c\rt)\det(Dh^{\ep}(z)) dz\nn\\
&\leq& \ep \int_{h^{\ep}(\Pi)} \det\lt(\lt[Du(y)\rt]_c\rt) dy\nn\\
&\leq& c\ep. 
\end{eqnarray}
\red{And thus we have established (\ref{obbeq51}).}

Similarly since $D\phi_u^{\ep}(z)=Du(h^{\ep}(z)) Dh^{\ep}(z)$ so 
\begin{eqnarray}
\label{gheq53}
\int_{\Pi} \det(D\phi_u^{\ep}(z)) dz&=&\int_{\Pi} \det\lt(Du(h^{\ep}(z))\rt) \det\lt(D h^{\ep}(z)\rt) dz\nn\\
&\leq& \int_{\Omega} \det\lt(Du(y)\rt) dy\nn\\
&\leq& \int_{\Omega} \det\lt(\lt[Du(y)\rt]_c\rt) dy\nn\\
&\leq& c.
\end{eqnarray}
Note by $\det(D\phi_u^{\ep}(z))
\overset{(\ref{gheq50})}{=}\det\lt(\lt[D\phi_u^{\ep}(z)\rt]_c\rt)+
\det\lt(\lt[D\phi_u^{\ep}(z)\rt]_a\rt)$. So 
\begin{eqnarray}
\label{mmeq4}
\int_{\Pi} \det\lt(\lt[D\phi_u^{\ep}(z)\rt]_c\rt) dz&\leq& 
\int_{\Pi} \lt|\det\lt(D\phi_u^{\ep}(z)\rt)\rt| dz+\int_{\Pi} \lt|\det\lt(\lt[D\phi_u^{\ep}(z)\rt]_a\rt)\rt| dz\nn\\
&\overset{(\ref{gheq52}),(\ref{gheq53})}{\leq}&c+c\ep\leq c.
\end{eqnarray}
Now $\|D\phi_u^{\ep}(z)\|\overset{(\ref{brq1})}{\leq} \sqrt{\det(\lt[D\phi_u^{\ep}(z)\rt]_c)}+ \sqrt{\lt|\det(\lt[D\phi_u^{\ep}(z)\rt]_a)\rt|}$ so 
\begin{eqnarray}
\|D\phi_u^{\ep}\|_{L^2(\Pi)}&\leq& c\lt(\int_{\Pi} \det(\lt[D\phi_u^{\ep}\rt]_c) dz \rt)^{\frac{1}{2}}
+c\lt(\int_{\Pi} \lt| \det\lt(\lt[D\phi_u^{\ep}\rt]_a\rt)\rt|\rt)^{\frac{1}{2}}\nn\\
&\overset{(\ref{gheq52}),(\ref{mmeq4})}{\leq}& c\text{ for all }\ep>0\nn
\end{eqnarray}
which establishes (\ref{obbeq53}).

So as previously explained letting $\phi_u$ be the weak limit of $\phi_u^{\ep_k}$ (for some subsequence $\ep_k$) we have 
that $\phi_u$ is holomorphic by Weyl's lemma. And by uniform convergence of $h^{\ep_k}$ to $h$ we have that 
$\phi_u=u\circ h$. In exactly the same way $\phi_v$ is holomorphic.\nl

\em Step 3. \rm We will show $h$ is a homeomorphism. 

\em Proof of Step 3. \rm From Step 2 we know $\phi_u=u\circ h$ is a holomorphic function. Recall $h$ is the 
uniform limit of monotone maps $h^{\ep}$. So as before by Kuratowski, Lacher and Whyburn \cite{mcau} we have that $h$ is monotone. So we know that for any $a\in w(\Omega)$, $h^{-1}(a)$ is 
connected. Suppose we can find two distinct points $p_1,p_2\in h^{-1}(a)$. Then as $\phi_u$ must be constant $h^{-1}(a)$ and this is an infinite set so $\phi_u$ must be constant on $w(\Omega)=\Omega'$. Since $h(\Omega')=\Omega$ and $\phi_u=u\circ h$ 
function $\phi_u$ can only be constant if $u$ is constant on $\Omega$ which is a 
contradiction. Hence $h^{-1}(a)$ consists of a single point for every $a\in \Omega$ and so $h$ is a homeomorphism.\nl 
 
\em Proof of Lemma completed. \rm Since $h$ is a homeomorphism we can define $w=h^{-1}$ and we then have the decompositions of $u$ and $v$ given in 
(\ref{obbeq30}) simply from the definition of $\phi_u$ and $\phi_v$. Also note that $w^{-1}=h\in W^{1,2}$. $\Box$

%
%

\begin{a1}
\label{LL6}
Suppose function $u\in W^{1,1}(\Omega:\R^2)$ has a decomposition of the form $u=\phi\circ w$ where 
$\phi$ is holomorphic then let $z_1,z_2,\dots \in w(\Omega)$ be the zeros 
of $\phi'$ and let $y_i=w^{-1}(z_i)$ we will show that for any open set $\OI\subset \subset \Omega\backslash \lt(\bigcup_{i=1}^{\infty} y_i\rt)$ we have $w\in W^{1,1}(\OI:\R^2)$.
\end{a1}

\em Proof of Lemma \ref{LL6}. \rm Since $\phi$ is a holomorphic the set of zeros of $\phi$ forms a countable set 
$z_1,z_2,\dots \in w(\Omega)$. So for any $z\in w(\Omega)\backslash \bigcup_{i=1}^{\infty} z_i$ there exists $r_z>0$ such that 
$\phi$ is a homeomorphism on $B_{r_z}(z)$. Let $y_i=w^{-1}(z_i)$ and let $y\in \Omega\backslash \bigcup_{i=1}^{\infty} y_i$. By continuity 
of $w$ on $\Omega$ there exists $\alpha_y>0$ such that $w(B_{\alpha_y}(y))\subset B_{r_{w(y)}}(w(y))$. Now 
$u(B_{\alpha_y}(y))=\phi(w(B_{\alpha_y}(y)))\subset \phi(B_{r_{w(y)}}(w(y)))$ and as $\phi$ is invertible on 
$\phi(B_{r_{w(y)}}(w(y)))$ we have that $w=\phi^{-1}\circ u$ on $B_{\alpha_y}(y)$. 
Since $\phi^{-1}$ is Lipschitz this implies $w\in W^{1,1}(B_{\alpha_y}(y))$. As 
any $\OI\subset \subset \Omega\backslash \lt(\bigcup_{i=1}^{\infty} y_i\rt)$ can be covered by a finite number of balls on which $w\in W^{1,1}$ so the result follows. $\Box$

%
%

\begin{a1}
\label{LL5}
Suppose $\Omega$ is a connected open set and $u\in W^{1,1}(\Omega:\R^2)$ is a mapping of 
integrable dilatation for which we have the decomposition $u=\phi\circ w$ 
where $\phi$ is a holomorphic and $w$ is a homeomorphism with \red{$w^{-1}\in W^{1,1}(w(\Omega):\R^2)$}, then $w\in W^{1,1}(\Omega:\R^2)$. 
\end{a1}

\em Proof of Lemma \ref{LL5}. \rm \red{By Lemma \ref{LL6} for any subset $\OI\subset \subset \Omega\backslash \lt(\bigcup_{i=1}^{\infty} w^{-1}(z_i)\rt)$ we have $u\in W^{1,1}(\OI)$. So by the chain rule $Du(x)=D\phi(w(x))Dw(x)$ for 
a.e.\ $x\in \OI$. Thus by Lemma \ref{LL1}
\begin{equation}
\label{gheq5}
\mu_{Dw(x)}=\mu_{Du(x)}\text{ for }a.e.\ x\in \OI.\nn
\end{equation}  
So $w$ is a homeomorphism of finite dilatation on $\OI$. Thus by Theorem 3.3 \cite{hen3}, $w^{-1}$ is a mapping of finite dilatation. By Theorem 1.2 \cite{hen3} $w$ satisfies 
\begin{equation}
\label{menq1}
\int_{\OI} \lt|Dw\rt| dx=\int_{w(\OI)} \lt|Dw^{-1}\rt| dx\leq \int_{u(\Omega)} \lt|D w^{-1}\rt| dx.
\end{equation} 
This holds for every $\OI\subset \subset \Omega$. Now we can define a vector field\footnote{Strictly speaking an 
equivalence class of vector fields that agree a.e.\ } 
$\varpi:\Omega\rightarrow M^{2\times 2}$ by defining $\varpi=Dw$ on $\OI$ where $\OI$ is any open set with $\OI\subset 
\subset \Omega$. So let $\OI_k$ be a sequence of sets $\OI_{k}\subset \OI_{k+1}$ with $\bigcup_{k} \OI_k=\Omega\backslash 
\lt(\bigcup_{i=1}^{\infty} w^{-1}(z_i)\rt)$. Now by (\ref{menq1}) we have 
$$
\int_{\Omega\backslash \lt(\bigcup_{i=1}^{\infty} w^{-1}(z_i)\rt)} \lt|\varpi\rt| dx=\lim_{k\rightarrow \infty} 
\int_{\OI_k} \lt|\varpi\rt| dx\leq \int_{u(\Omega)} \lt|D u^{-1}\rt| dx.
$$
Thus $\varpi\in L^1(\Omega)$, note also for any 
$\phi\in C_c( \Omega\backslash \lt(\bigcup_{i=1}^{\infty} w^{-1}(z_i)\rt) :\R^2)$ since $\mathrm{Spt} \phi\subset \OI_k$ for large enough $k$ so 
\begin{equation}
\int w \mathrm{div} \phi \;dz=\int_{\OI_k} w \mathrm{div} \phi \;dz
=-\int_{\OI_k} Dw\cdot \phi\; dz= -\int \varpi\cdot \phi \;dz.
\end{equation}
Thus $\varpi$ does indeed serve as the distributional derivative of $w$ and hence 
$w\in W^{1,1}(\Omega\backslash \lt(\bigcup_{i=1}^{\infty} w^{-1}(z_i)\rt):\R^2)$, which is 
easily seen to imply $w\in W^{1,1}(\Omega:\R^2)$.} $\Box$\nl

%
%

\subsection{Proof of Theorem completed.} By Lemma \ref{LL4} there exists homeomorphism $w:\Omega\rightarrow \R^2$ and holomorphic 
function $\phi_u, \phi_v$ such that $u=\phi_u\circ w$ and $v=\phi_v\circ w$. By Lemma \ref{LL5} we have that $w\in W^{1,1}(\Omega:\R^2)$. 
\red{So by the chain rule
\begin{eqnarray}
Dv(z)&=&\PPI\lt(\phi_v'\lt(w\lt(z\rt)\rt)\rt)Dw(z)\nn\\
&=&\PPI\lt(\phi_v'\lt(w\lt(z\rt)\rt)\rt)\lt(\PPI\lt(\phi_u'\lt(w\lt(z\rt)\rt)\rt)\rt)^{-1} \PPI\lt(\phi_u'\lt(w\lt(z\rt)\rt)\rt) Dw(z)\nn\\
&=&\PPI\lt(\phi_v'\lt(w\lt(z\rt)\rt)\rt)\lt(\PPI\lt(\phi_u'\lt(w\lt(z\rt)\rt)\rt)\rt)^{-1} Du(z)\nn\\
&=& \PPI\lt(\frac{\phi_v'\lt(w\lt(z\rt)\rt)}{ \phi_u'\lt(w\lt(z\rt)\rt)}\rt) Du(z).\nn
\end{eqnarray}
So} let $\psi(z):=\frac{\phi_v'(z)}{\phi_u'(z)}$, $\psi$ is a meromorphic function and $Dv(z)=\PPI(\psi(w(z)))Du(z)$ and 
hence we have established (\ref{gheq99}). 

Let $\lt\{z_1,z_2,\dots \rt\}$ be the set of zeros of $\phi_u'$, then by \red{Stoilow} decomposition the 
branch set of $u$ is given by $B_u=w^{-1}\lt(\lt\{z_1,z_2,\dots \rt\}\rt)$. So as $\lt\{z_1,z_2,\dots \rt\}$ 
are also the set of poles of $\psi$ we have that they are contained in $w(B_u)$. $\Box$\nl

%
%

\section{Proof of Corollaries.}

\subsection{Proof of Corollary \ref{CC2}} By Theorem \ref{T2} there exists meromorphic $\psi$ and homeomorphism 
$w\in W^{1,1}$ such that $Du(z)=\PPI(\psi(w(z)))Dv(z)$. Now since $S(Du(z))=S(Dv(z))$ for a.e.\ $z$ so $\det(Du(z))=\det(P(\psi(w(z))))\det(Dv(z))$ and thus 
\begin{equation}
\label{jjeq1}
\lt|\psi(w(z))\rt|=1\text{ for }a.e.\ z\in \Omega. 
\end{equation}
Let $\lt\{z_1,z_2,\dots \rt\}$ be the set of singularities of $\psi$ on $w(\Omega)$. Note that since $\lt\{z_1,z_2,\dots \rt\}$ 
only have cluster points on $\partial w(\Omega)$ so $\Omega\backslash \lt\{w^{-1}(z_1),w^{-1}(z_2),\dots \rt\}$ is open. 
So for any $\varpi\in \Omega\backslash \lt\{w^{-1}(z_1),w^{-1}(z_2),\dots \rt\}$ we have $r>0$ such that 
$B_r(\varpi)\subset \Omega\backslash \lt\{w^{-1}(z_1),w^{-1}(z_2),\dots  \rt\}$. Hence as $w$ is an open mapping and 
$\psi$ is an open mapping on $w(B_{r}(\varpi))$ we have that $\psi\circ w$ is open on $B_r(\varpi)$. Thus by 
(\ref{jjeq1}) $\psi\circ w$ is constant $B_r(\varpi)$. So we have $\zeta_0\in \mathbb{C}\cap\lt\{z:\lt|z\rt|=1\rt\}$ 
such that $\psi(w(z))=\zeta_0$ for any $z\in B_r(\varpi)$. Let 
$$
\Lambda:=\lt\{z\in \Omega\backslash \lt\{w^{-1}(z_1),w^{-1}(z_2),\dots  \rt\}: \psi(w(z))=\zeta_0 \rt\}.
$$
By continuity of $\psi\circ w$ the set $\Lambda$ is closed. By the argument above $\Lambda$ is also open and so 
as $\Omega\backslash \lt\{w^{-1}(z_1),w^{-1}(z_2),\dots  \rt\}$ is a connected open set and we know $\Lambda=\Omega\backslash \lt\{w^{-1}(z_1),w^{-1}(z_2),\dots  \rt\}$. Thus $Du(z)=\PPI(\zeta_0)Dv(z)$ on $\Omega$ and hence we have established (\ref{bvbb50}). $\Box$

\subsection{Proof of Corollary \ref{CC4}} 

Let $z_1,z_2,\dots \in w(\Omega)$ be the zeros of $\phi_u'$. Recall $B_u$ is the branch set of $u$ and note 
$$
\lt\{w^{-1}(z_1),w^{-1}(z_2),\dots \rt\}=B_u
$$ 
and if $\zeta\in \Omega\backslash \lt\{w^{-1}(z_1),w^{-1}(z_2),\dots \rt\}$ then by continuity of $u$ 
there exists $r_{\zeta}>0$ 
such that $\phi_u^{-1}$ is well defined on $u(B_{r_{\zeta}}(\zeta))$. Hence 
\begin{equation}
w(z)=\phi_u^{-1}(u(z))\text{ for any }z\in B_{r_{\zeta}}(\zeta).\nn
\end{equation}
Now by Theorem \ref{T2} the set of poles of $\psi$ is contained in $w(B_u)$. So for every 
$\zeta\in \Omega\backslash B_u$ there exists $\lambda_{\zeta}>0$ such that $\psi$ is holomorphic on $w(B_{\lm_{\zeta}}(\zeta))$. Let $\tau_{\zeta}=\min\lt\{r_{\zeta},\lm_{\zeta}\rt\}$. Then we have that 
\begin{equation}
\label{gv2}
Dv(z)=\PPI(\psi(\phi_u^{-1}(u(z))))Du(z)\text{ for }z\in B_{\tau_{\zeta}}(\zeta).
\end{equation}
Now $\psi\circ \phi_u^{-1}$ is analytic on $u(B_{\tau_{\zeta}}(\zeta))$ and since $\zeta\not \in B_u$ assume $\tau_{\zeta}$ is 
chosen small enough to that $u_{\lfloor B_{\tau_{\zeta}}}$ is homeomorphic and so $u(B_{\tau_{\zeta}}(\zeta))$ is simply connected. 
So by Cauchy's theorem there exists holomorphic $\xi$ defined on $u(B_{\tau_z}(z))$ such that $\xi'(z)=\psi(\phi_u^{-1}(z))$ for all $z\in u(B_{\tau_{\zeta}}(\zeta))$. 
Thus from (\ref{gv2}) we have that 
 \begin{equation}
Dv(z)=\PPI(\xi'(u(z)))Du(z)\text{ for }z\in B_{\tau_{\zeta}}(\zeta).\nn
\end{equation}
Now let $p=\mathrm{Re}(\xi)$, $q=\mathrm{Im}(\xi)$ so $\xi'(x,y)=p_x(x,y)+iq_x(x,y)$. Thus  
$$
\PPI\lt(\xi'(u(z))\rt)=\lt(\begin{matrix} p_x(u(z)) & -q_x(u(z))\\ q_x(u(z)) & p_x(u(z))) \end{matrix}\rt)=\lt(\begin{matrix} p_x(u(z)) & p_y(u(z))\\ q_x(u(z)) & q_y(u(z))) \end{matrix}\rt).
$$
Let $\varphi(x,y)=(p(x,y),q(x,y))$. $D\varphi(x,y)=\lt(\begin{smallmatrix} p_x & p_y \\ q_x & q_y \end{smallmatrix}\rt)$ so 
$Dv(z)=D\varphi(u(z))Du(z)$ for all $z\in B_{\tau_{\zeta}}(\zeta)$ which establishes the corollary. $\Box$

%
%

\subsection{Proof of Corollary \ref{CC3}}

\red{We can find a collection of bounded simply connected sets $\UI_k$ such that $\UI_k\subset \UI_{k+1}$ for each $k$ and 
$\bigcup_{k=1}^{\infty} \UI_k=\Omega$. For example if $\zeta_0\in \Omega$ we can take $\UI_k$ to be the 
connected component of $\Omega\cap B_{2^{k}}(\zeta_0)$ containing $\zeta_0$. It is easy to see 
$\Omega\cap B_{2^k}(\zeta_0)$ is simply connected because for any $\sigma\not \in \Omega\cap B_{2^k}(\zeta_0)$ the function $\frac{1}{z-\sigma}$ has a primitive on $\Omega\cap B_{2^k}(\zeta_0)$.} 
 
\red{Given complex function $u$ of a complex variable we can, as before define 
$$
\ti{u}(x,y)=\lt(\mathrm{Re}(u(x+iy)),\mathrm{Im}(u(x+iy))\rt).
$$
As we have noted $\lt[D \ti{u}(x,y)\rt]_C=\lt[\frac{\partial u}{\partial z}(x+iy)\rt]_M$ where $\lt[a+ib\rt]_M=\lt(\begin{smallmatrix}  a & -b \\ b & a\end{smallmatrix}\rt)$. Thus for 
any fixed $k$, the function $\ti{u}$ on $\UI_k$ has the hypothesis to apply Corollary 2. Since $u$ is a homeomorphism $B_u=0$ and thus by Corollary \ref{CC4} for any $\zeta\in \UI_k$ there exists $r_{\zeta}>0$ such that 
$v=\phi^{\zeta}\circ u$ on $B_{r_{\zeta}}(\zeta)$. Also because $u$ is a homeomorphism we know $u(\UI_k)$ is simply connected. And 
for any $\zeta\in \Omega$, $\phi^{\zeta}$ is defined on $u(B_{r_{\zeta}}(\zeta))$. Pick $y_0\in u(\Omega)$ for any other 
$y_1\in u(\UI_k)$ we can find an injective path $\Gamma$ starting at $y_0$ and ending at $y_1$. So 
$\phi^{u^{-1}(y_0)}$ can be analytically continued along the path to $y_1$. Thus by the Mondromy theorem there exists analytic 
function $\psi_k$ defined on all of $u(\UI_k)$ such that $\psi_k=\phi^{u^{-1}(y)}\text{ on }u(B_{r_{u^{-1}(y)}}(u^{-1}(y)))\text{ for any }y\in u(\UI_k)$. 
Thus 
$$
v(z)=\phi^{u^{-1}(y)}(u(z))=\psi_k(u(z))\text{ on }z\in B_{r_{u^{-1}(y)}}(u^{-1}(y))\text{ for any }y\in u(\UI_k). 
$$
Hence $v=\psi_k\circ u$ on $\UI_k$. Now if $l>k$ arguing in the same way we have the decomposition $v=\psi_l\circ u$ on $\Pi_l$ for analytic function 
$\psi_l$ defined on $u(\Pi_l)$. However since $u(\Pi_l)$ is simply connected by the Mondromy theorem $\psi_l$ is the extension 
of $\psi_k$ to $u(\Pi_l)$. So we can define analytic function $\psi$ on $u(\Omega)$ by letting 
$\psi(z)=\psi_k(z)$ when $z\in u(\Pi_k)$. Then $v=\psi\circ u$ on $u(\Omega)$. $\Box$}

\section{Counter example}
\label{example}

The example below that shows the sharpness of the hypotheses of Theorem \ref{T2} and also shows the sharpness of the two dimensional case of Theorem 1 of \cite{lor20}. It first appeared in 
more general form in Section 5 of \cite{lor20}. Subsequent to \cite{lor20} being 
accepted we learned of the article \cite{cia} where an example is given (attributed to Herv\'{e} Le Dret and 
one of that article's referees) that showed that (\ref{bvbb51.5}) does not imply (\ref{bvbb50}) 
for arbitrary functions in $W^{1,2}$. The example presented in \cite{cia} does not show the sharpness of space 
of functions of integrable dilatation however it is in spirit not unrelated to the example of 
\cite{lor20}. As the example of \cite{lor20} is simple to 
describe we present the two dimensional version of it for the convenience of the reader.

\em Example 1. \rm Let $Q_1:=\lt\{z\in \R^2:\lt|z\rt|_{\infty}<1\rt\}$. Define 
\begin{equation}
u(x_1,x_2):=\lt\{\begin{array}{ll} (x_1,x_2 x_1)&\text{ for } x_1>0\nn\\
(x_1,-x_2 x_1)&\text{ for } x_1\leq 0\end{array}\rt.
\end{equation}
and for some $\theta\in (0,2\pi)$
\begin{equation}
v(x_1,x_2):=\lt\{\begin{array}{ll} (x_1\cos \theta-x_1 x_2\sin \theta,x_1\sin\theta+x_1x_2\cos\theta)&\text{ for } x_1>0\nn\\
(x_1,-x_2 x_1)&\text{ for } x_1\leq 0\end{array}\rt.
\end{equation}
\red{Note $u$, $v$ are Lipschitz}. Now that for $x_1\leq 0$
\begin{equation}
\label{fte1}
Du(x)=Dv(x)=\lt(\begin{matrix} 1&0\\
-x_2&-x_1\end{matrix}\rt).
\end{equation}
And for $x_1>0$
\begin{eqnarray}
\label{fte2}
\na v(x)&=&\lt(\begin{matrix} \cos\theta-x_2\sin\theta&-x_1\sin\theta\\
\sin\theta+x_2\cos\theta&x_1\cos\theta\end{matrix} \rt)
=\lt(\begin{matrix} \cos\theta&-\sin\theta\\
\sin\theta&\cos\theta\end{matrix}\rt)\lt(\begin{matrix} 1& 0\\
x_2&x_1
\end{matrix}\rt),
\end{eqnarray}
\begin{eqnarray}
\label{fte2.3}
\na u(x)=\lt(\begin{matrix} 1& 0\\
x_2&x_1
\end{matrix}\rt).
\end{eqnarray}
So from (\ref{fte1}), (\ref{fte2}), (\ref{fte2.3}) it is clear there is no $R$ such that $\na v(x)=R \na u(x)$ for $x\in Q_1$. Now note that $\det(\na u(x))=x_1$ for all $x\in Q$ and $\lt|\na u(x)\rt|^2
=\lt(1+x_2^2+x_1^2\rt)^{\frac{1}{2}}$ so defining $L(x):=\lt|\na u(x)\rt|^2/\det(\na u(x))=x_1^{-1}\lt(1+x_2^2+x_1^2\rt)^{\frac{1}{2}}$. So it is clear that $\int_{Q_1} L(z) dz=\infty$ and 
thus it follows that there can be no continuous relation between $\na u$ and $\na v$ on a dense \red{connected} open subset of $\Omega$ without the assumption of 
integrability of the dilatation and thus Theorem \ref{T2} is optimal. $\Box$

\section{Appendix}

\begin{a1}
\label{LL10}
\red{Let $\Omega$ be an open domain, let $u\in W^{1,1}(\Omega:\R^2)$ be a homeomorphism,  then }
\begin{equation}
\red{\det(Du)\in L^1_{loc}(\Omega).}\nn
\end{equation}
\end{a1}
\em Proof of Lemma \ref{LL10}. \rm \red{Since $u$ is approximately differentiable a.e. in $\Omega$ it satisfies hypothesis (a) of Theorem 
1 of \cite{haj}. Now for any set $S\subset \subset \Omega$ we can find a bounded open set $\Pi\subset\subset \Omega$ with 
$S\subset \Pi$.} 

\red{Now by Theorem 2 \cite{haj} we can redefine $u$ on a set of zero measure $N\subset \Omega$ to create a new function $\ti{u}$ that satisfies 
\begin{equation}
\label{fgg1}
\int_{E} g\circ \ti{u}\lt|\det(D \tilde{u})\rt| dx=\int_{\R^2} g(y) N_{\ti{u}}(y,E) dy\text{ for any }E\subset \Omega
\end{equation} 
where $N_{\ti{u}}(y,E)=\mathrm{card}(\ti{u}^{-1}(y),E)$ and $g$ is any measurable function. So let $\Pi'=\Pi\backslash N$ and $g=\chara_{u(\Pi')}$. So 
$N_{\ti{u}}(y,\Pi')=N_{u}(y,\Pi')=1$ for any $y\in u(\Pi')$. Note as $\Pi$ is a bounded and $u$ is homeomorphism on $\Omega$ so 
$u(\Pi)$ is bounded. Hence} 
$$
\red{\int_{\Pi} \lt|\det(Du)\rt| dx=\int_{\Pi'} \lt|\det(D \tilde{u})\rt| dx\overset{(\ref{fgg1})}{=}\int_{\R^2} \chara_{u(\Pi')} dy=\lt|u(\Pi')\rt|<\infty}. \;\;\;\;\Box
$$

\end{document}